\documentclass[]{elsarticle}
\usepackage{moreverb}

\usepackage{hyperref}
\usepackage{times}
\usepackage{calc}
\usepackage{graphics}
\usepackage{epsfig}
\usepackage{amssymb}
\usepackage{amsmath}
\usepackage{flafter}
\usepackage{float}
\usepackage{latexsym}
\usepackage{psfrag}
\usepackage{rotating}
\usepackage{longtable}
\usepackage{subfigure}
\usepackage{color}

\renewcommand{\vec}[1]{\bold{#1}}

\newcommand{\Vnabla}{\vec{\nabla}}

\begin{document}
\begin{frontmatter}
\journal{arxiv}
\title{Arbitrary Order Solutions for the Eikonal Equation using a Discontinuous Galerkin Method}

\author[1]{David Flad\corref{cor1}} 
\author[2]{Aniruddhe Pradhan} 
\author[1]{Scott Murman}
\address[1]{NASA Ames Research Center, Moffet Field, CA, USA}
\address[2]{Department of Mechanical Engineering, University of Michigan, Ann Arbor, MI 48109, USA}
\cortext[cor1]{Corresponding author, Universities Space Research Association, david.g.flad@nasa.gov}

\begin{abstract}

We provide a method to compute the entropy-satisfying weak solution to the eikonal equation in an arbitrary-order polynomial space.  The method uses an artificial viscosity approach and is demonstrated for the signed distance function, where exact solutions are available. The method is designed specifically for an existing high-order discontinuous-Galerkin framework, which uses standard convection, diffusion, and source terms. We show design order of accuracy and good behavior for both shocks and rarefaction type solutions. Finally the distance function around a complex multi-element airfoil is computed using a high-order-accurate representation.
\end{abstract}


\end{frontmatter}

\section{Introduction}\label{sec:intro}
The eikonal equation $|\Vnabla s| = 1/f$, is the Hamilton-Jacobi equation of classical mechanics and optics.
Solutions to the eikonal equation are relevant in a wide range of applications such as interface reinitialization \cite{Sussman:94,Peng:99,sukumar2001,wang2003,peirce2008,Peirce2015,fechter2015,Nguyen2019}, mesh generation \cite{Bagnerini:01,Wang:05,Persson:05,Zhao2016}, image segmentation \cite{Sofou:03,Alvino:07} or first-arrival travel times \cite{Mo:02,Taillandier:09,Lan:13}. 
Many fast efficient algorithms were developed to solve the equation such as the finite-difference Fast Marching \cite{Tsitsiklis:95, Sethian:96}, Fast Sweeping \cite{zhao:04,detrixhe:13}, and Group Marching \cite{Kim:01} methods, along with finite-element Discontinuous-Galerkin methods \cite{Zhang:11}.
The most efficient methods achieve $\mathcal{O}(N)$ complexity. In order to achieve this efficiency, the algorithms rely on the hyperbolic nature of the equation, resulting in an exploitable causality for the solution update. For high-order methods, structured non-oscillatory algorithms are typically used.

These schemes for the eikonal equation require specialized algorithms, which are often difficult to implement efficiently in a parallel environment, and do not easily generalize to unstructured, high-order methods.

Following the philosophy of Tucker \cite{Tucker:03}, a discontinuous Galerkin solver for the eikonal equation was developed in \cite{murman2015} for body-fitted, general unstructured meshes achieving high-order accuracy. The motivation for this approach is that, while some algorithmic efficiency is lost, the eikonal equation is solved in parallel, leveraging the same data structures and algorithms as the physics simulation for which the solution is required. High-order, globally smooth solution are often beneficial when needed within the solution of the primary simulation, such as in the implementation of turbulence models. This greatly simplifies its development and usage, as standard methods for solving partial differential equations based upon convection, diffusion, and source terms are utilized.

As in \cite{murman2015}, this work first focuses on computing an approximate signed distance function $s$, resulting as $f=1$ and $s=0$ on boundaries $\Gamma$. We seek to leverage our work on discontinuous-Galerkin spectral-element methods (DGSEM) for the Navier-Stokes equations \cite{diosady2013,diosady2014,diosady2015} to develop an eikonal solver. In order to align the equation with the solver strategy, we solve a modified system of equations for a decoupled system of gradients. This reduces the stiffness of the system, but does not provide convergence to the entropy solution. Thus we couple the original equation weakly through a source term. Finally, we introduce an artificial viscosity that decreases with mesh resolution, providing the appropriate weak solution in the presence of shocks.

Due to an efficient three-dimensional tensor-product implementation, the cost for a residual evaluation of the method can be made independent of order of accuracy through the use of software optimization\cite{diosady2013}. The framework runs efficiently on modern high-performance computing clusters, is of arbitrary order of accuracy, and can handle general unstructured meshes.

\section{Modified eikonal system of equations}\label{sec:eikoanl}
The eikonal equation is given by,
\begin{align}\label{eikonal}
|\Vnabla s|&= \frac{1}{f} \\
s&=g \quad on  \ \ \Gamma \nonumber.
\end{align}
For a constant force term $f=1$ and $g=0$, the solution $s$ is the shortest distance to the boundary $\Gamma$ for any given point. 
It is informative to view this equivalently as the minimal time needed to travel from a given point to the boundary $\Gamma$. Then we realize that the gradient may be seen as the velocity inverse $\Vnabla s = 1/\vec{v}$ and $f$ is the velocity magnitude at any given point. Thus the eikonal equation is likewise a solution to the problem of a curve moving into its own normal direction with speed $f$. One interesting application is to solve the eikonal equation with varied velocity functions, for example moving under the influence of its own curvature, as is often the case in multi-phase problems. An extensive discussion is found in the book of Sethian \cite{sethianbook}.

\begin{figure}
    \centering
    \includegraphics[width=0.8\textwidth]{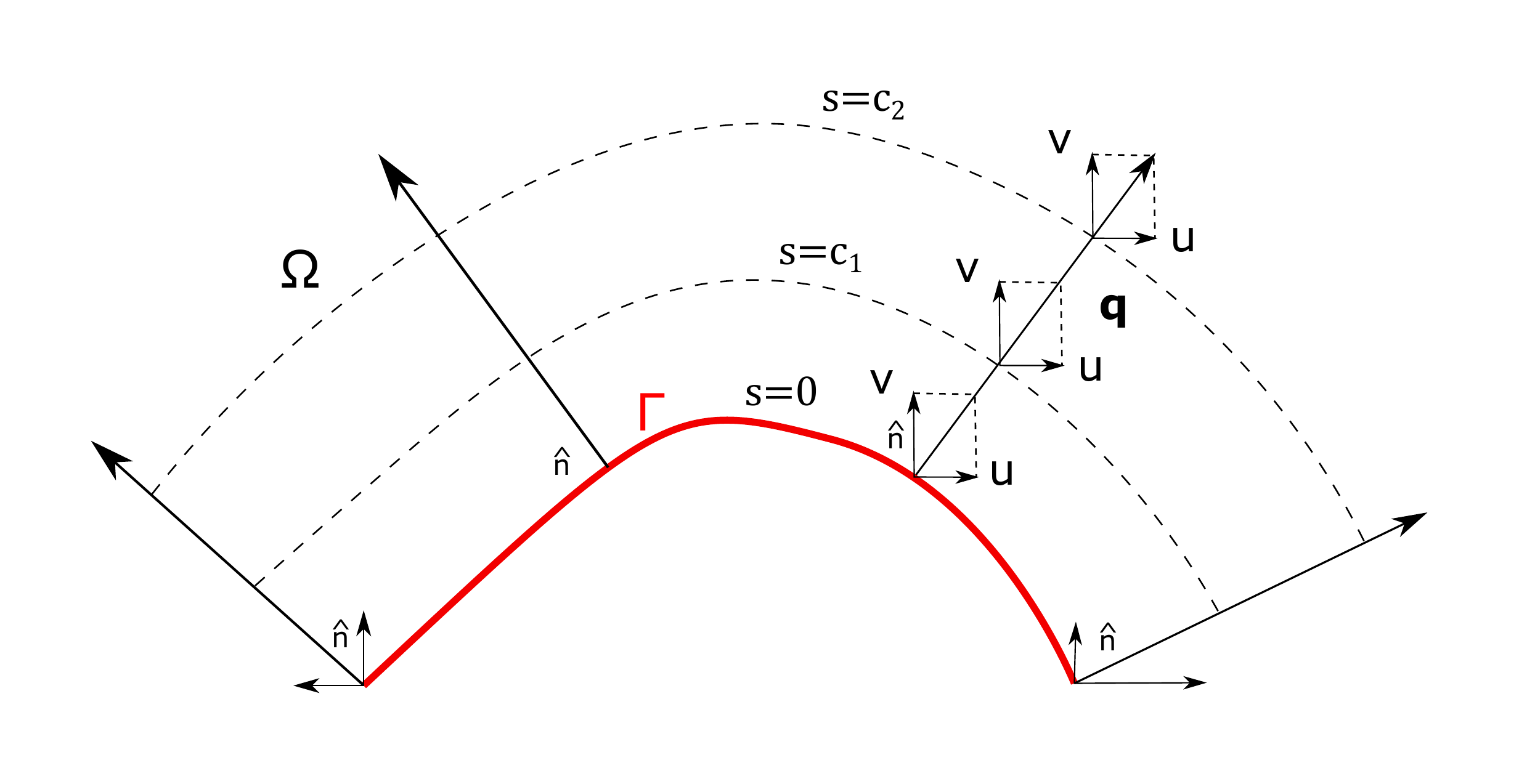}
    \caption{ Gradient vector components $u,v,w$ advect along the characteristic direction.}
    \label{fig:eikonal_intro}
\end{figure}

\begin{figure}
    \centering
    \includegraphics[width=0.8\textwidth]{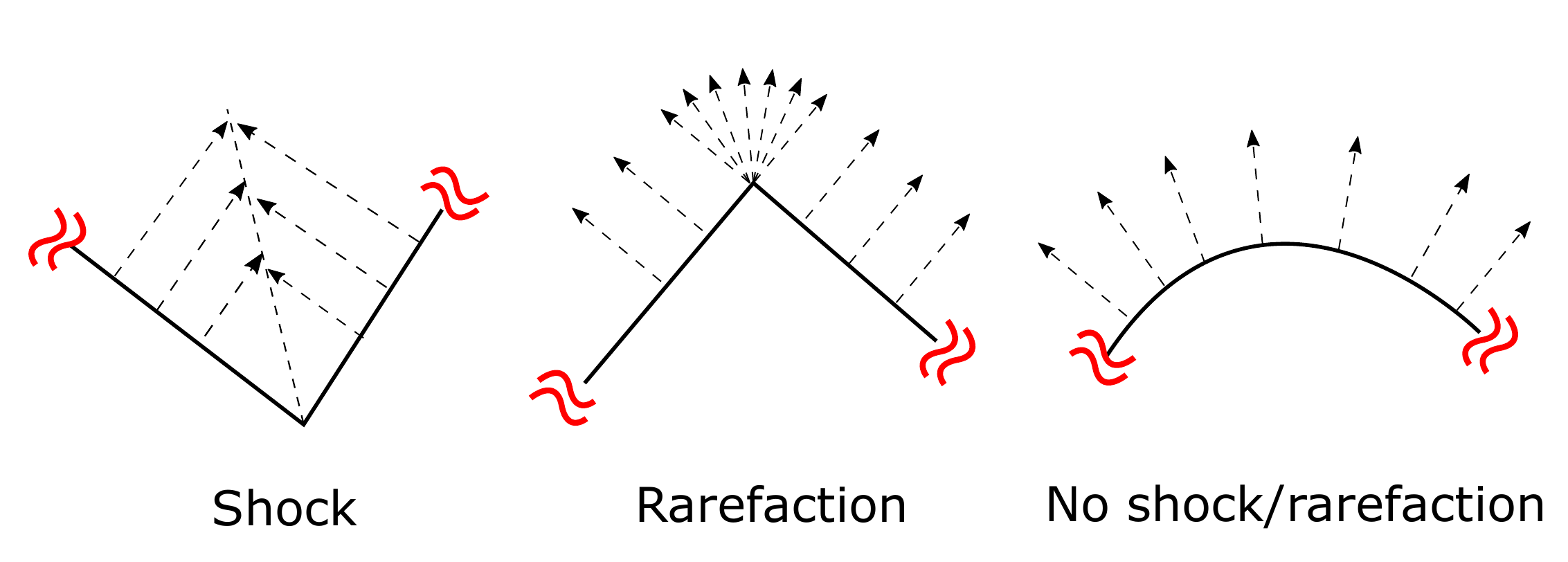}
    \caption{Topology of the surface can lead to shock formation, rarefaction and no shock/rarefaction cases.}
    \label{fig:eikonal_scenarios}
\end{figure}

To cast the equation into a form that is better suited to a discontinuous-Galerkin framework, we square the eikonal equation
\begin{align}\label{eikonalSquared}
\Vnabla s \cdot \Vnabla s&= \frac{1}{f^2} \\
s&=g \quad on  \ \ \Gamma \nonumber.
\end{align}
As a consequence the solution is limited to one side of the boundary, {\emph{i.e.}},  we limit ourselves to computing positive or negative values for the distance $s$. Now we introduce the gradient vector $\vec{q} = [u,v,w]^T = \Vnabla s$. As the original equation was found to be very difficult to converge to a steady state solution using a high-order solver \cite{murman2015,Hartmann2014}, we decouple the solution of the gradient field and the scalar distance function. To get a set of partial differential equations for the gradients, we take the derivative of equation \eqref{eikonalSquared} with respect to Cartesian coordinate directions $\vec{x} = [x,y,z]^T$:
\begin{equation}
\Vnabla (\vec{q}\cdot\vec{q}) = \Vnabla (\frac{1}{f^2}).
\end{equation}
Exemplary for the x-direction,
\begin{align}
(uu)_x + (vv)_x + (ww)_x &= \left(\frac{1}{f^2}\right)_x \\
\textnormal{with} \quad v_x &= s_{yx}=s_{xy}=u_y
\end{align}
follows
\begin{align}
uu_x + vu_y + wu_z &= \left(\frac{1}{2f^2}\right)_x.
\end{align}
Our system to solve for the gradients is then found as
\begin{align}\label{eqn:gradSys}
(\vec{q} \cdot \Vnabla)\vec{q}  = \Vnabla ({1\over 2f^2}).
\end{align}
This system is physically intuitive, showing that the gradient vector components $u,v,w$ are transported along the direction of $\vec{q}$, as depicted in figure \ref{fig:eikonal_intro}.
Now we continue by rewriting \eqref{eqn:gradSys} into a conservative form with source terms as
\begin{align}\label{burgers}
\Vnabla \cdot \vec{F} &= \Vnabla \frac{1}{2f^2} + \vec{q} \ \Vnabla \cdot \vec{q}\\ \nonumber
\vec{q} \cdot \vec{n} &= 1 \quad on \ \ \Gamma,
\end{align}
where $\vec{n}$ is the boundary normal vector and $\vec{F}^l$ direction-wise flux vector
\begin{align}
\vec{F}^1 = \begin{pmatrix}
uu\\
uv\\
uw
\end{pmatrix}
\vec{F}^2 = \begin{pmatrix}
vu\\
vv\\
vw
\end{pmatrix}
\vec{F}^3 = \begin{pmatrix}
wu\\
wv\\
ww
\end{pmatrix},
\end{align}
i.e. a three-dimensional form of Burgers equation. As shown in figure \ref{fig:eikonal_scenarios}, the topology of the domain can lead to the formation of shocks and rarefactions. For cases that do not contain shocks or strong rarefactions, solving equation \eqref{burgers} and subsequently integrating is sufficient. 

Instead of a system of three differential equations, we may also integrate the original equation \eqref{eikonalSquared} cast in divergence form with substituted gradients by $\vec{q}$
\begin{align}\label{eikonalSqDiv}
\Vnabla \cdot (s\vec{q}) = \frac{1}{f^2} + s \Vnabla \cdot \vec{q},
\end{align}
the form used in \cite{murman2015}. However, for cases where the solution is not smooth, leading to shocks forming in the gradients equation, the solution becomes multi-valued and depends on the initial solution. To avoid this we want to find a weak, entropy satisfying solution for the original equation \eqref{eikonalSquared}, as discussed in \cite{sethianbook}. Hence, we have to couple the system \eqref{burgers} with equation \eqref{eikonalSqDiv}, which forces the gradients to be those of the scalar valued variable $s$. 

Before introducing a coupling method, the resulting system of differential equations is modified by adding a viscosity term, that we can use later to enforce convergence to the weak solution
\begin{align}
\Vnabla \cdot (s\vec{q}) &= \frac{1}{f^2} + s \Vnabla \cdot \vec{q} +  \mu \Vnabla \cdot \Vnabla s \label{sysPart1} \\ 
\Vnabla \cdot \vec{F}&= \Vnabla \frac{1}{2f^2} + \vec{q} \ \Vnabla \cdot \vec{q} + \mu \Vnabla \cdot \Vnabla \vec{q}\label{sysPart2}
\end{align}

There exist multiple ways we might couple \eqref{sysPart1} into \eqref{sysPart2} by re-substituting one of the $\vec{q}$ terms with $\Vnabla s$. However, none was found to be robust, making convergence of high-order approximations difficult. Thus we use a weak coupling, introducing a source term on the right-hand side, which enforces the condition $\Vnabla s = \vec{q}$. 

Following the discussion above, that $s$ is an arrival time, we can non-dimensionalize the equations with a reference velocity $U_{ref}$ and length $L_{ref}$, leading to a non-dimensional number $\mathcal{P} = \frac{L_{ref}}{\mu U_{ref}}$. The final non-dimensional\footnote{For clarity we do not introduce special notation for non-dimensional variables} system we propose is
\begin{align}
\Vnabla \cdot (s\vec{q}) &= \frac{1}{f^2} + s \Vnabla \cdot \vec{q} +  \frac{1}{\mathcal{P}} \Vnabla \cdot \Vnabla s\nonumber \\ 
\Vnabla \cdot \vec{F}&= \Vnabla \frac{1}{2f^2} + \vec{q} \ \Vnabla \cdot \vec{q} + \frac{1}{\mathcal{P}} \Vnabla \cdot \Vnabla \vec{q}  - \lambda(\vec{q} - \Vnabla s), \label{sysEikonal}
\end{align}
where $\lambda=\frac{1}{L_{ref}f}$ is a multiplier that ensures dimensional consistency. In equations \eqref{sysEikonal}, we note that when shifted to the left hand side, the coupling source term $\lambda\Vnabla s$, is analogues to a pressure term, balanced by the source term $\lambda q$.
For the remainder of the paper, we look for a solution of the distance function where $f=1$ and $\Vnabla f =0$, to demonstrate the method and evaluate the accuracy against known solutions.
We define the viscosity as a function of space and solution, $\mu=\mu(s,\vec{q})$. The specific form we found to work well is similar to what was used in \cite{Tucker:03, Hartmann2014}
\begin{align}\label{eqn:artvisc}
\mu = \frac{c \Delta\sqrt{\bar{s}}}{\sqrt{L_{ref}f}},
\end{align}
where $\Delta=Vol^{1/3}/N$ (element volume divided by discretization order) is a reference length scale specific to the discretization, ensuring the viscosity goes to zero as $\Delta \downarrow 0$, $c$ is a user defined constant, and $\bar{s}$ is an element-wise integral mean of the distance function. In this way the artificial viscosity is reduced near the body where $\bar{s} \downarrow 0$. While we find equation \eqref{eqn:artvisc} efficient for scaling the amount of viscosity needed, a nonlinear shock detector is needed to regain the full advantage of the high-order discretization in smooth regions.  We examine the novel algorithm behavior with the simpler artificial viscosity formulation here, and future work will extend the approach from Murman et al. \cite{murman2017} to the eikonal system.

\section{Numerical method}
We discretize equation \eqref{sysEikonal} in the discontinuous-Galerkin framework described in \cite{diosady2014,diosady2015}. The solver allows for arbitrarily high order of accuracy. Using a tensor-product discretization and exploiting vectorization on modern high-performance computing systems, the cost of a residual evaluation is independent of the chosen order of accuracy \cite{diosady2013}. The computational domain is split into non-overlapping finite elements $\Omega$, with a Riemann solver utilized at discontinuous interfaces $\partial\Omega$. 

First, the system of equation \eqref{sysEikonal} is re-written in a compact flux-divergence form as
\begin{align}
    \Vnabla \cdot (\vec{G}-\vec{V})&= \vec{S},
    \label{sysEikonal_FD}
\end{align}
with the advective, viscous and source vectors given by
\begin{align}\label{eqn:eikonalSys}
\vec{G^l}= \begin{pmatrix}
            q^ls\\
            q^lu - g_1\delta_{1l} s\\
            q^lv- g_1\delta_{2l} s\\
            q^lw- g_1\delta_{3l} s
            \end{pmatrix} \quad
\vec{V^l}= \begin{pmatrix}
         \mu \Vnabla^l s\\
         \mu \Vnabla^l \vec{q}
         \end{pmatrix} \\ \nonumber
\vec{S}= \begin{pmatrix}
         \frac{1}{f^2}+s\Vnabla \cdot \vec{q} g_2 |\vec{q}|\\
         \Vnabla \frac{1}{2f^2} + \vec{q} \ \Vnabla \cdot \vec{q}%/(g_3|\vec{q}|) 
         - g_1\vec{q}
         \end{pmatrix} \quad l=1,2,3
\end{align} with $\delta$ the Kronecker symbol. With \eqref{eqn:eikonalSys}, two numerical parameters $g_1$ and $g_2$ are introduced, their individual purposes is discussed in section \ref{sec:switches}.

The test function and basis function are constructed using tensor-products from 1D Lagrange polynomials of degree $p$, resulting in a formally $p+1 = N^{\text{th}}$-order discretization, defining a polynomial space $\mathcal{L}_2$ for all elements $\Omega$. The discretized weak form then reads, find $[s,\vec{q}] \in \mathcal{L}_2$
\begin{align}
<(\mathcal{G^*}+\mathcal{V^*})\cdot \vec{n},\varphi>_{\partial\Omega}
-<{\mathcal{G}}+{\mathcal{V}},\Vnabla\varphi>_\Omega &= <{\mathcal{S}},\varphi>_\Omega \quad \forall \ \varphi \in \mathcal{L}_2.
\label{discEikonal}
\end{align}
where we denote by ${\mathcal{G,V,S}}$ the discrete flux and source terms. For the inter-element advective numerical flux, we use a Lax-Friedrich flux function,
\begin{equation}
\mathcal{G^*}=1/2(\mathcal{G^+}+\mathcal{G^-} - \gamma ([s^+,\vec{q}^+]-[s^-,\vec{q}^-])).
\end{equation}
As in \cite{murman2015}, we find using a global eigenvalue $\gamma = 1/f$ to improve convergence. For the viscous flux we use the second method of Bassi and Rebay \cite{bassi2}. To solve the nonlinear system, we use the diagonal regularization method described in \cite{ceze2016}, with a Newton-Krylov linear solver.

\subsection{Control parameters}\label{sec:switches}
In our numerical method, there exist two control parameters, denoted $g_1,g_2$ in \eqref{eqn:eikonalSys}. The purpose and behavior is: 
\begin{enumerate}
\item[i] $g_1 = \lambda$ whenever $\mu>0$.  This is usually the case in practical simulations as the viscosity is needed for any shape that has concave sections.   $g_1 = 0$ otherwise, so that the coupling can be switched off without loss of accuracy, reducing the overall stiffness of the system.
\item[ii] $g_2 = 1$ when $\Vnabla \cdot \vec{q} > 0$.  This is necessary in rarefaction cases. For negative  divergence of $\vec{q}$, the solution may include shocks and the solution we look for does not maintain $|q|=1$, thus $g_2=\frac{1}{|\vec{q}|}$.
\end{enumerate}
% , iii) $g_3$ is generally $0$ as $g_2$ is found sufficient for any case including viscosity ($g_1=1$, i.e. the system is coupled) and only used for the rarefaction test case of distance around a box, see section \ref{sec:box}.

\subsection{Dirichlet boundary condition}
On the boundary $\Gamma$, we prescribe $s=0$ and $\vec{q} = -\vec{n}$, with $\vec{n}$ the outward pointing normal vector of the discretized domain boundary. Thereby, in order to achieve full order of the scheme, the boundary representation is required to be at least $N+1$ order accurate, such that the prescribed boundary condition for the gradient vector $\vec{q}$ is of order $N$. This Dirichlet boundary condition is enforced weakly by setting
\begin{align}
    \mathcal{G^*} = \mathcal{G}([0, -\vec{n}]^T) \quad \text{and} \\
    \mathcal{V^*} = \mathcal{V^-} \quad \text{on} \ \ \Gamma,
\end{align}
as commonly done for discontinuous-Galerkin schemes.

\subsection{Artificial far-field boundary condition}
As we use a finite domain we need to set a boundary condition in the far-field $\Pi$, where no physical boundary is present. If we had no viscosity, the equation would be purely hyperbolic with only one eigenvalue. Besides the viscous, terms we find it sufficient for the proposed system to use the interior state
\begin{align}
    \mathcal{G^*} = \mathcal{G^-} \quad \text{and} \\
    \mathcal{V^*} = \mathcal{V^-} \quad \text{on} \ \ \Pi.
\end{align}

\section{Numerical experiments}

The numerical experiments in this chapter are set up to demonstrate formal accuracy of the scheme of up to order $N=16$ (cylinder), handling of concave shapes resulting in shock-like solutions (sinusoidal), two convex shapes to demonstrate the ability to do rarefaction-like solutions (square, NACA0012 airfoil) and a complex geometry to demonstrate a combination of the above (30P30N multi-element-airfoil). Plane geometry is chosen for the merit of better judgement of the quality of the results. Our meshing strategy uses an iso-parametric mapping of a reference domain with high-order splines. As discussed above, geometry representation is chosen to be one order higher than the solution $N+1$, except noted otherwise. Reference length and velocity are chosen as $1$ in all test cases.

\subsection{Distance from a cylinder}
\begin{figure}
    \centering
    \includegraphics[width=0.9\textwidth, trim=0 5 0 30 ,clip]{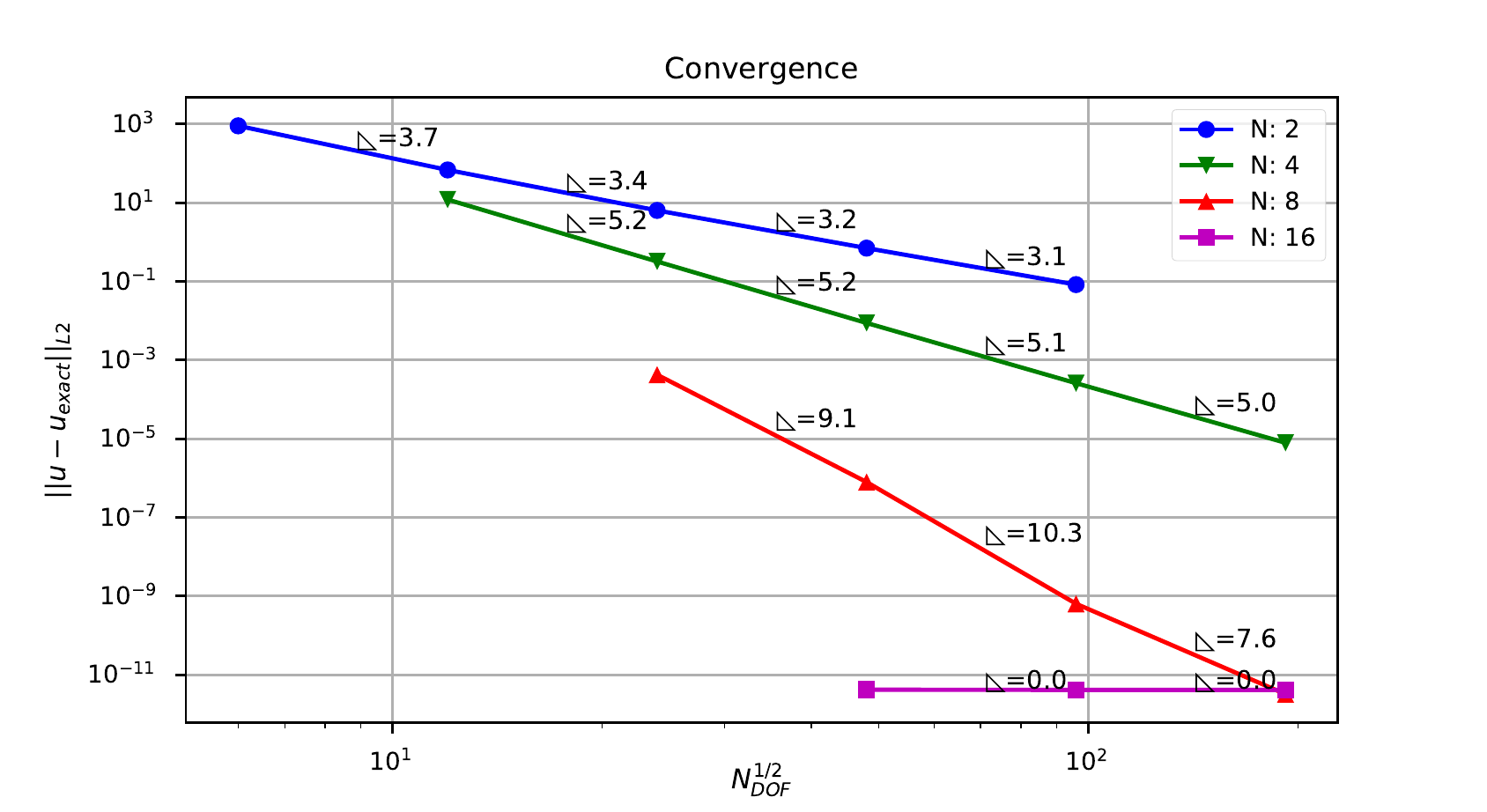}\\
        \includegraphics[width=0.9\textwidth, trim=0 5 0 30 ,clip]{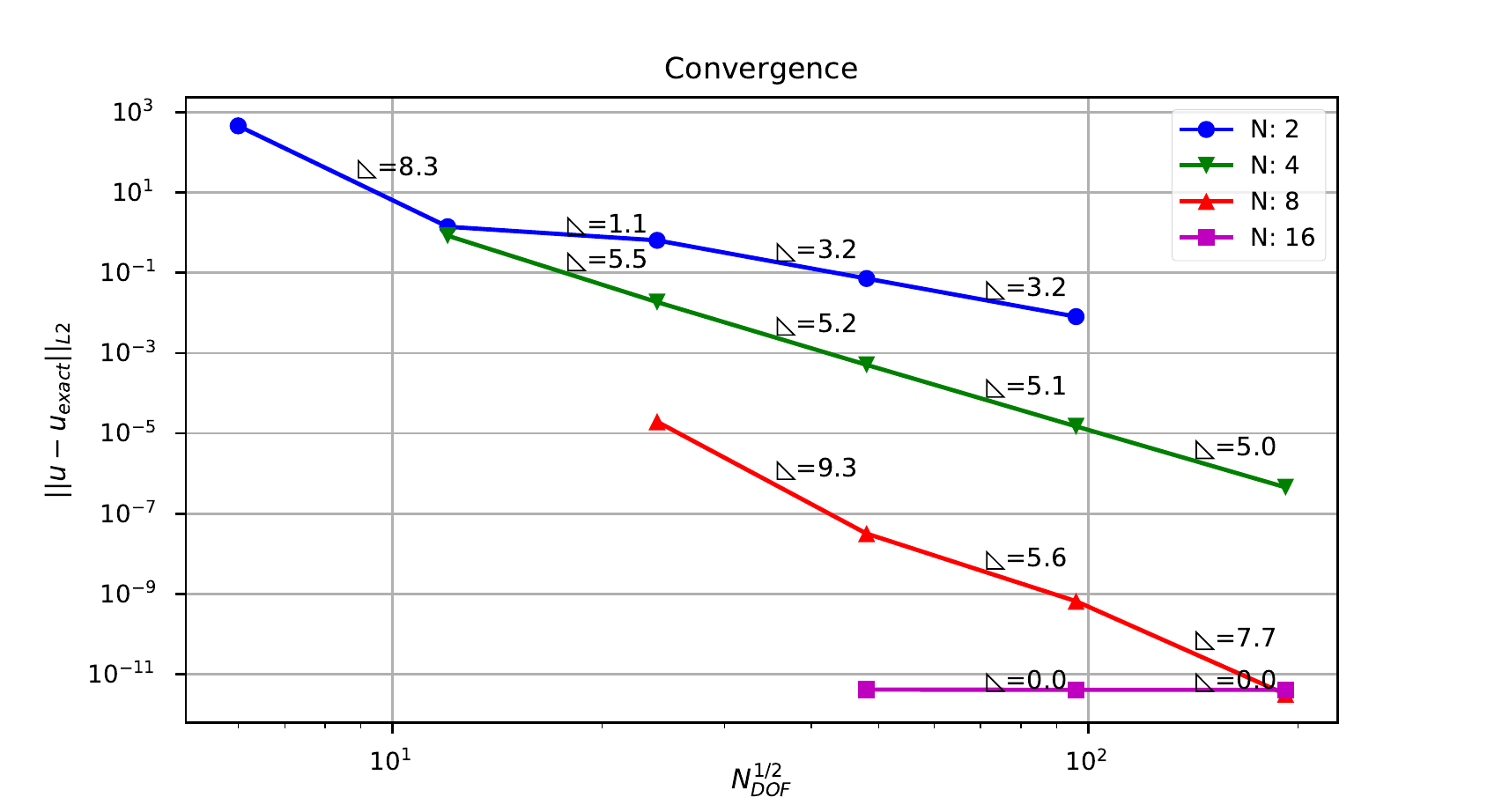}
    \caption{$L^2$ error of distance function around a cylinder. Top: without scaling of curvature ($g_2=0$), bottom: proposed method, with scaling of curvature}
        \label{fig:convtest}
\end{figure}

A cylinder is used to demonstrate formal order of accuracy of the code. For this purpose, the artificial viscosity is set to $\mu=0$. All results converge to the prescribed precision of $1\times e^{-10}$, without requiring any spectral prolongation from a lower order discretization. The error for the $16^{th}$-order approximation is already machine epsilon for the coarsest discretization of $3\times3$ elements, in radial and circumferential direction, mesh dimensions are doubled in both directions from here up to $48\times48$ elements. Figure \ref{fig:convtest} shows the $L^2$-norm error of a distance function around a cylinder for approximation orders $N=2,4,8,16$. For all orders, the convergence is found to be at least $N+1$, one order more than the discretization and in accordance with the boundary representation for $s$. As found in \cite{murman2015}, we see the benefit of using higher order approximations for this problem. In addition to the proposed method, the results with all control parameters off is plotted in figure \ref{fig:convtest}. It shows the clearer convergence behavior, but the error for each resolution and approximation order is higher.

\subsection{Distance from a sinusoidal}
The distance from a sinusoidal function at the bottom, and a plane wall at the top, is computed. This leads to a discontinuity at the lower channel half vertical centerline due to the concave shape of the sinusoidal, and a discontinuity at the horizontal center line with a triplet point in the center of the computational domain. An exact solution is not readily available for the problem, and we will asses the solutions only qualitatively here, with the purpose of demonstrating the shock-capturing capabilities of the proposed scheme. All solutions are converged to machine epsilon, without requiring spectral prolongation from lower order solutions. The artificial viscosity coefficient is calibrated using a 1D shock tube and set to $c=0.9$ for all cases. 

\begin{figure}
    \centering
    \includegraphics[height=8cm, trim=1000 200 1100 200 ,clip]{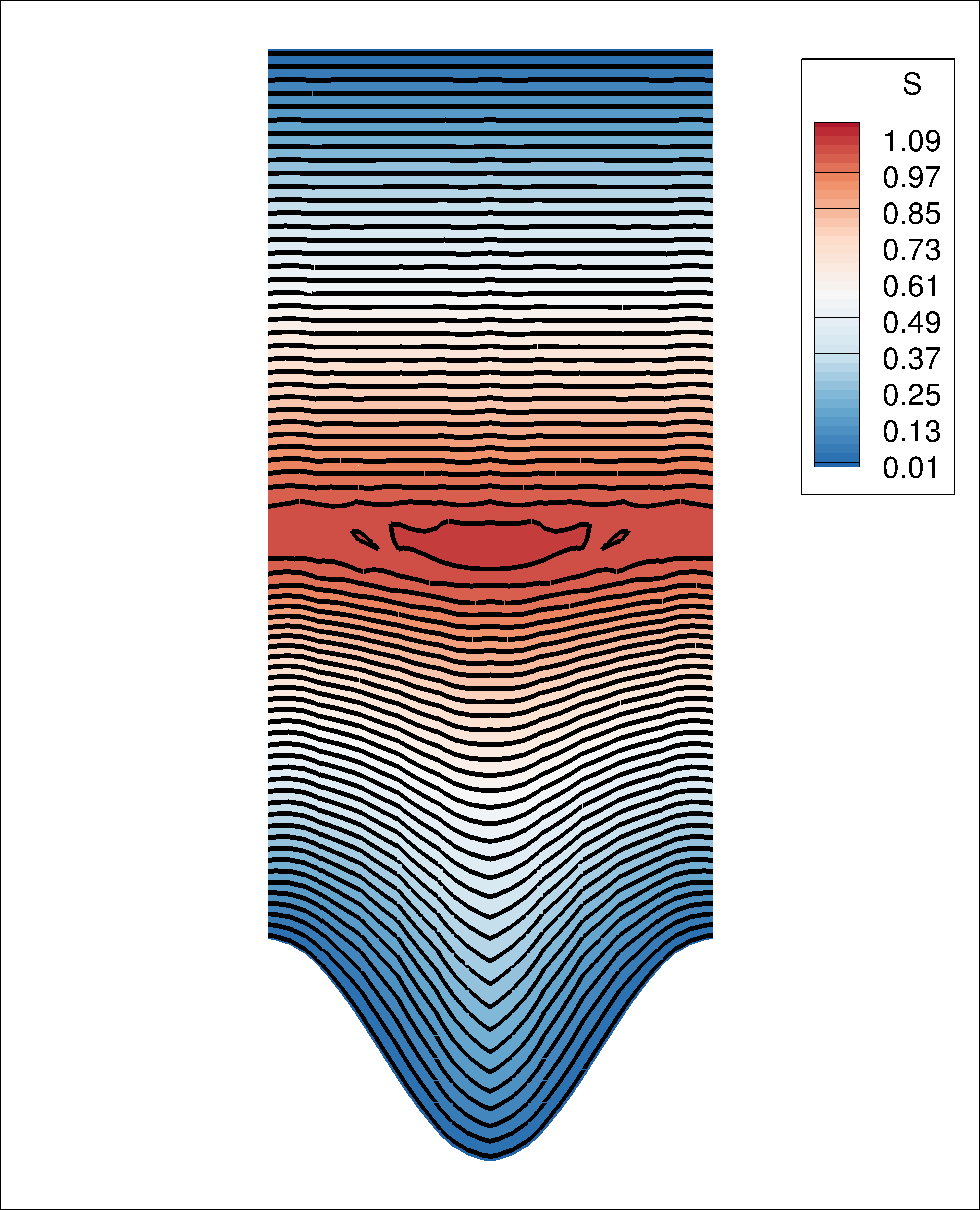}
    \includegraphics[height=8cm, trim=1000 200 1100 200 ,clip]{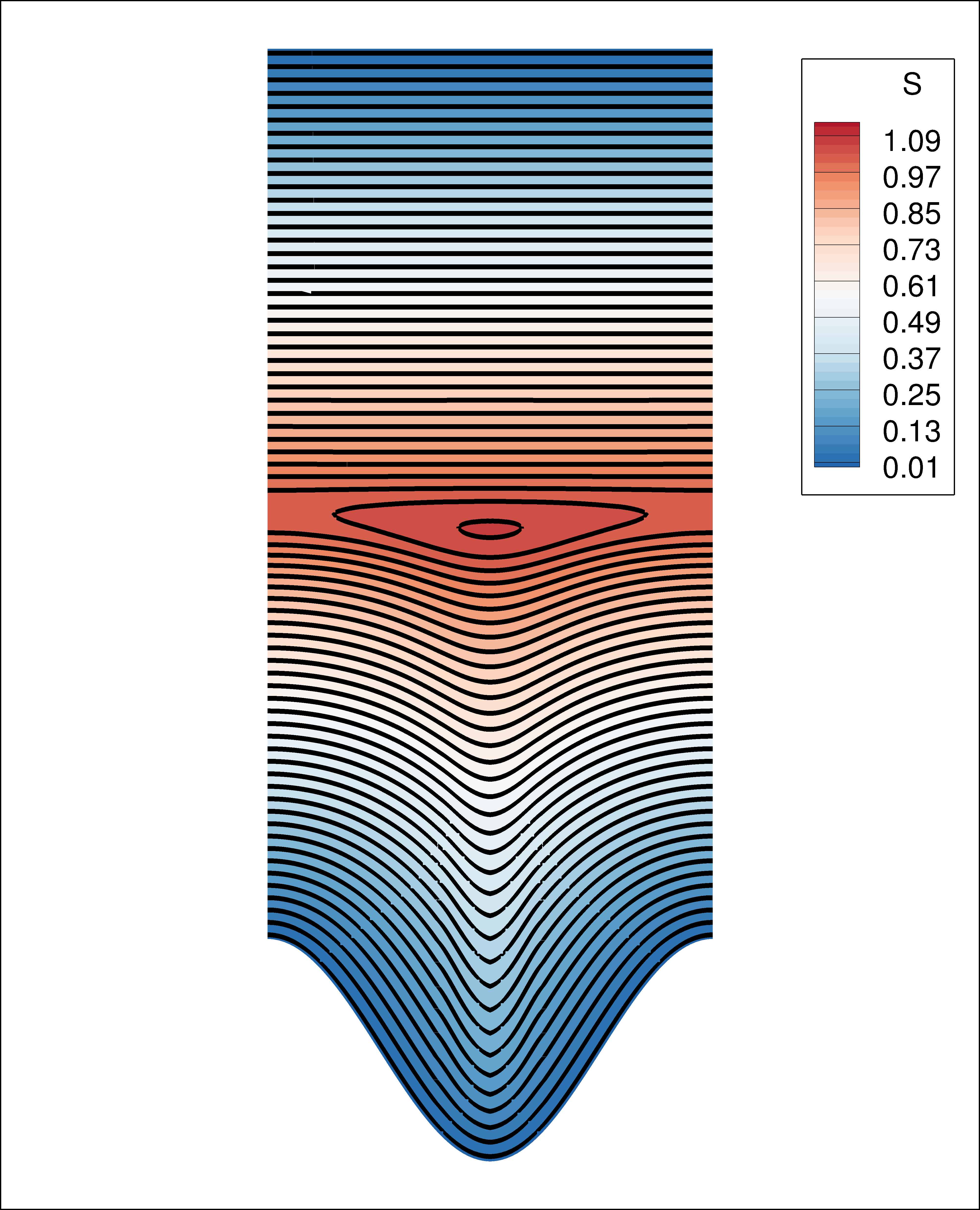}
    \includegraphics[height=8cm, trim=1000 200 100 200 ,clip]{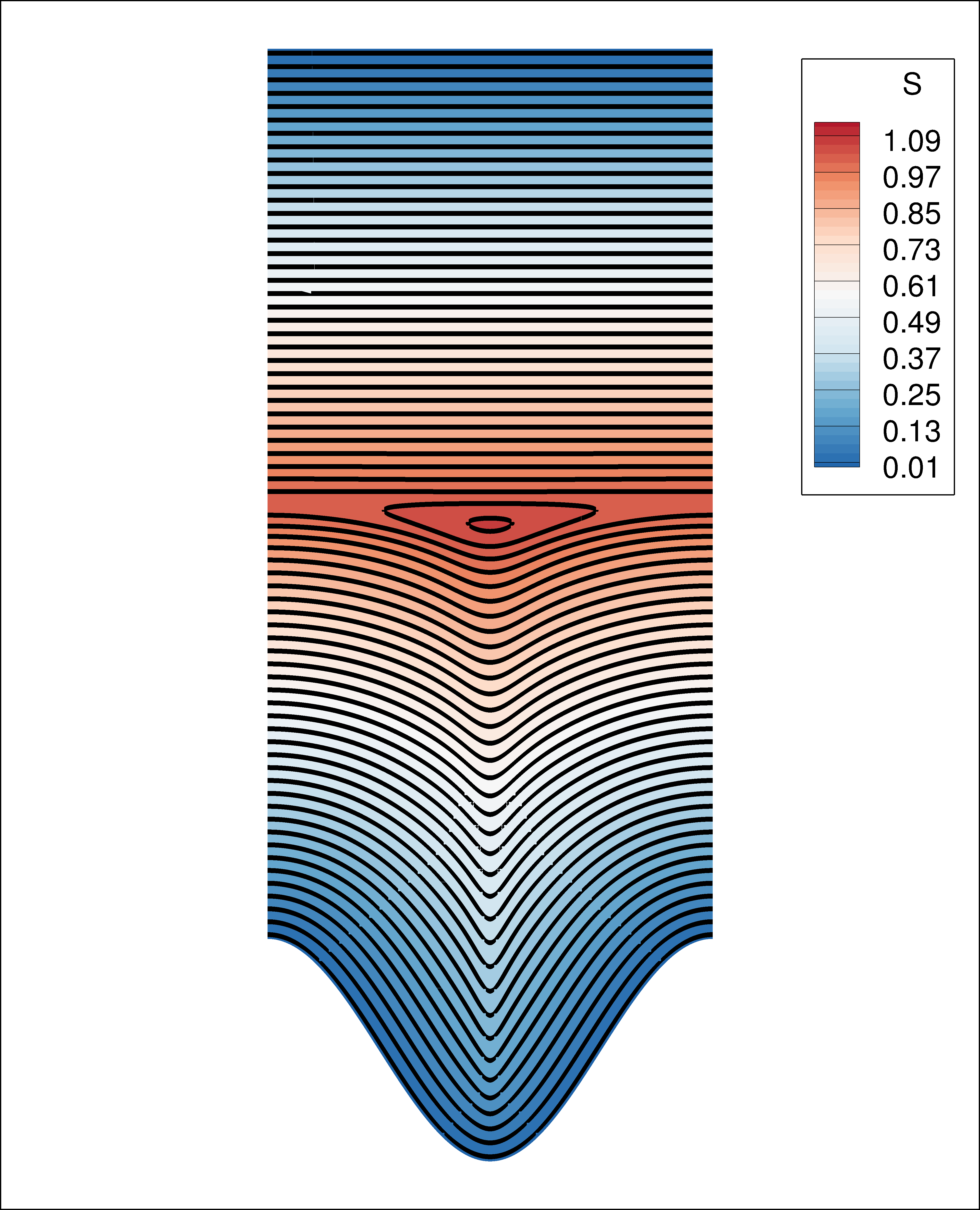}
    \caption{Solution for the sinusoidal channel problem for orders $N=2,4,8$ on $10\times30$ elements in horizontal and vertical direction, respectively.}
    \label{fig:sinusRefined}
\end{figure}

\begin{figure}
    \centering
        \includegraphics[height=8cm, trim=1000 200 1100 200 ,clip]{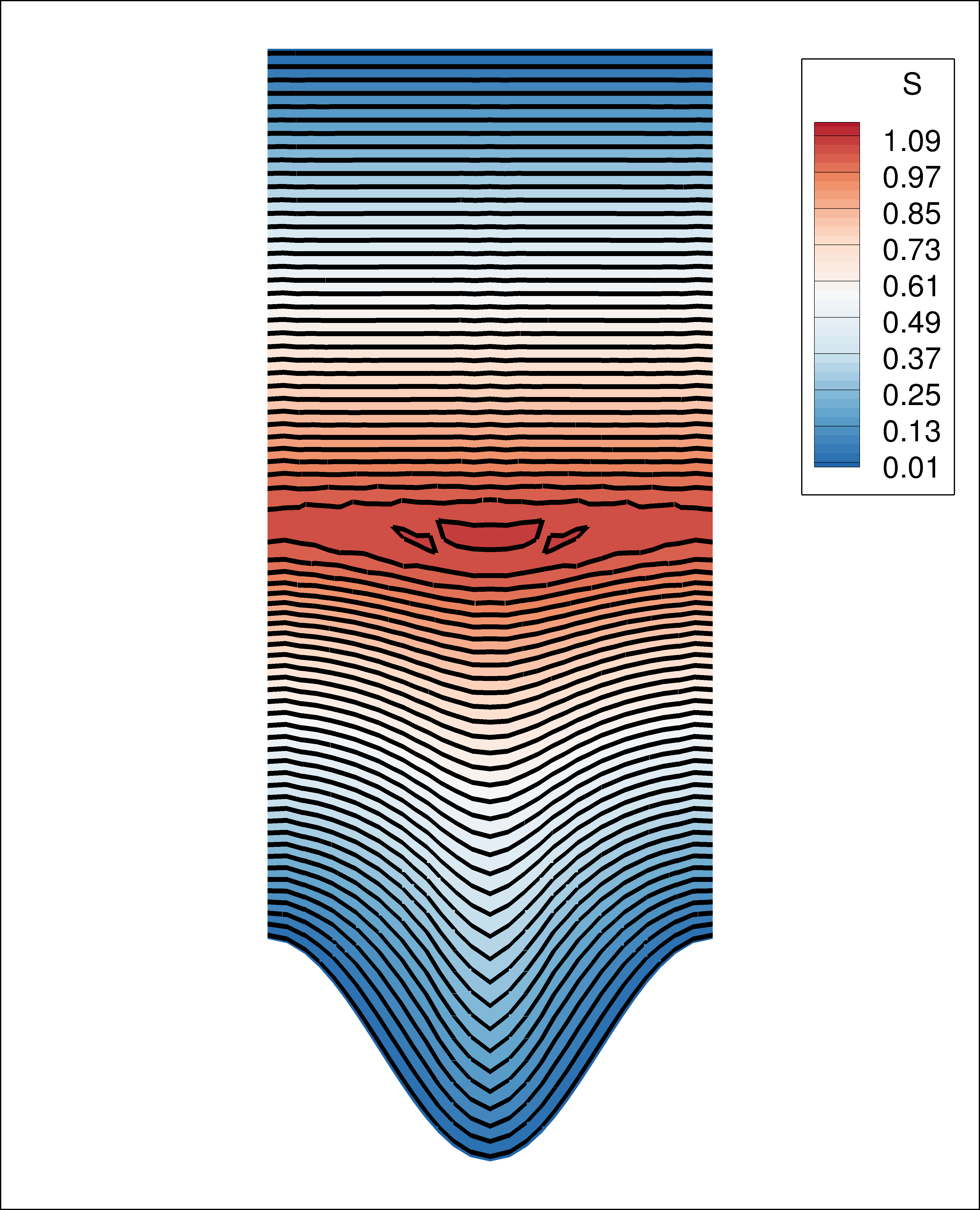}
        \includegraphics[height=8cm, trim=1000 200 1100 200 ,clip]{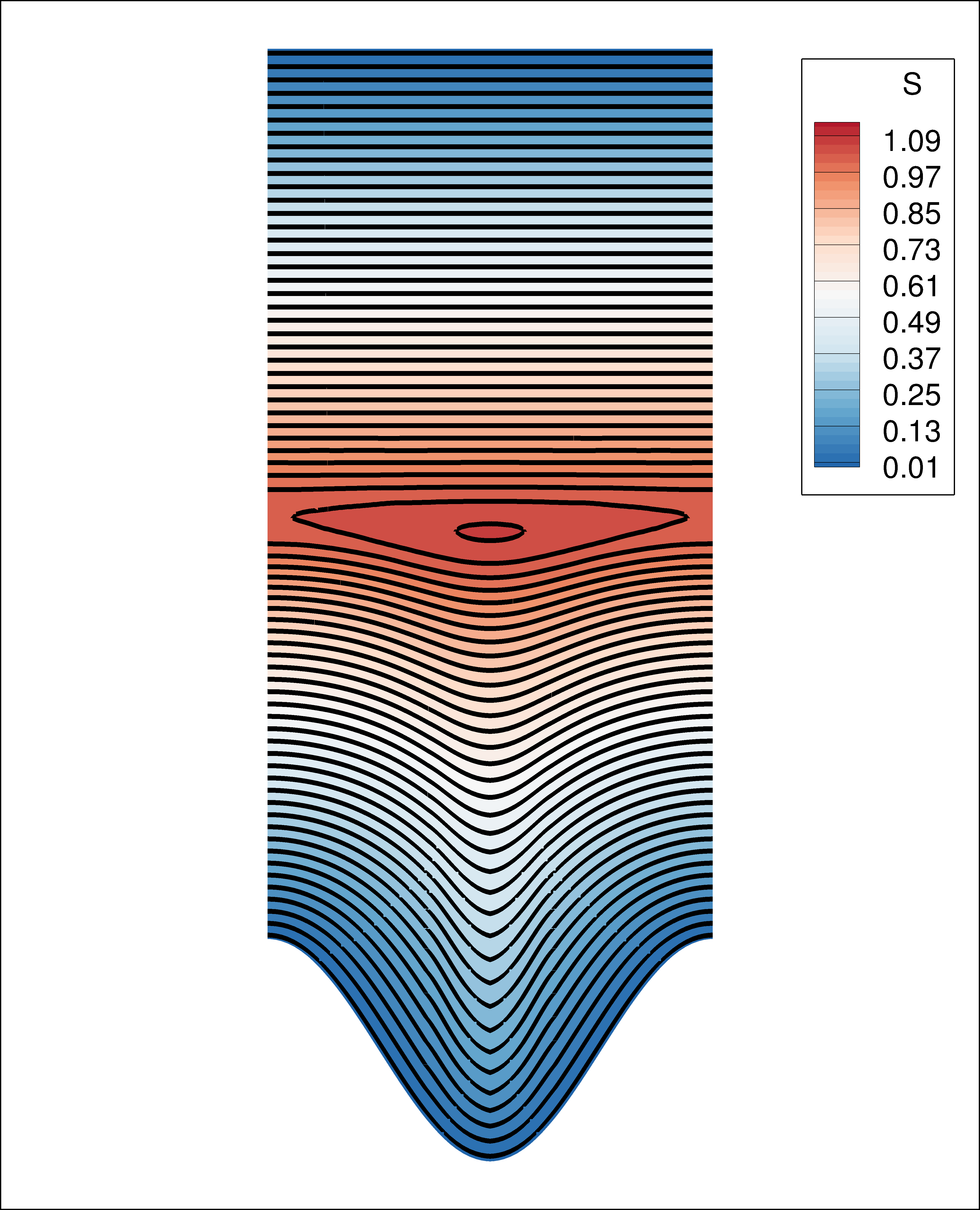}
    \includegraphics[height=8cm, trim=1000 200 100 200 ,clip]{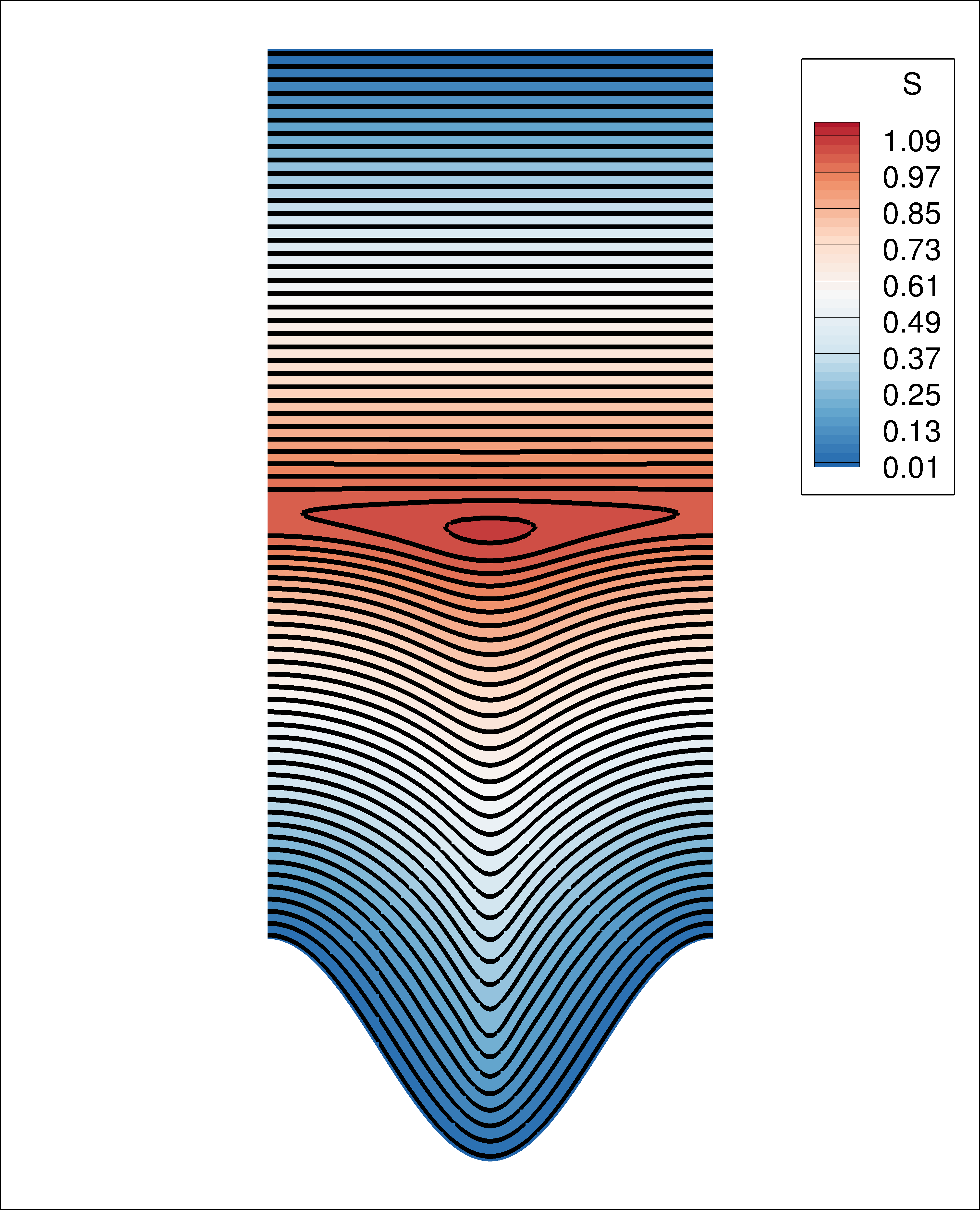}
    \caption{Solution for the sinusoidal channel problem for order $N=2,4,8$ with total DOF per direction $32\times64$, each.}
    \label{fig:sinusSameDOF}
\end{figure}

Figure \ref{fig:sinusRefined} shows the $N=2,4,8$ order solutions of the problem computed on a $10\times30$ element mesh in the horizontal and vertical directions, respectively. All solutions are found to converge to the entropy solution. The $N=2$ solution is significantly under-resolved, and does not capture well the vertical and horizontal discontinuity. With increasing resolution by means of approximation order, the discontinuities are captured with increasing detail. For order $N=8$, an almost sharp representation of the discontinuous solution in the lower channel half is found. The horizontal shock and triplet point are also captured better, relative to the lower order results.
Figure \ref{fig:sinusSameDOF} shows the results for orders $N=2,4,8$, while keeping the total number of degrees of freedom (DOF) constant at $32\times64$. We see that the high-order elements provide increased resolution relative to the same number of low-order DOF, capturing both discontinuities increasingly sharper.
In order to improve the scheme, $h-p$ adaptation, or similar, is required to properly account for discontinuities using high-order methods, which is beyond the scope of the current work. However, we still clearly see the benefit of high-order elements even in simulations with discontinuities. 

\subsection{Distance from a square}
\begin{figure}
    \centering
    \includegraphics[width=0.3\textwidth, trim=80 100 100 10  ,clip]{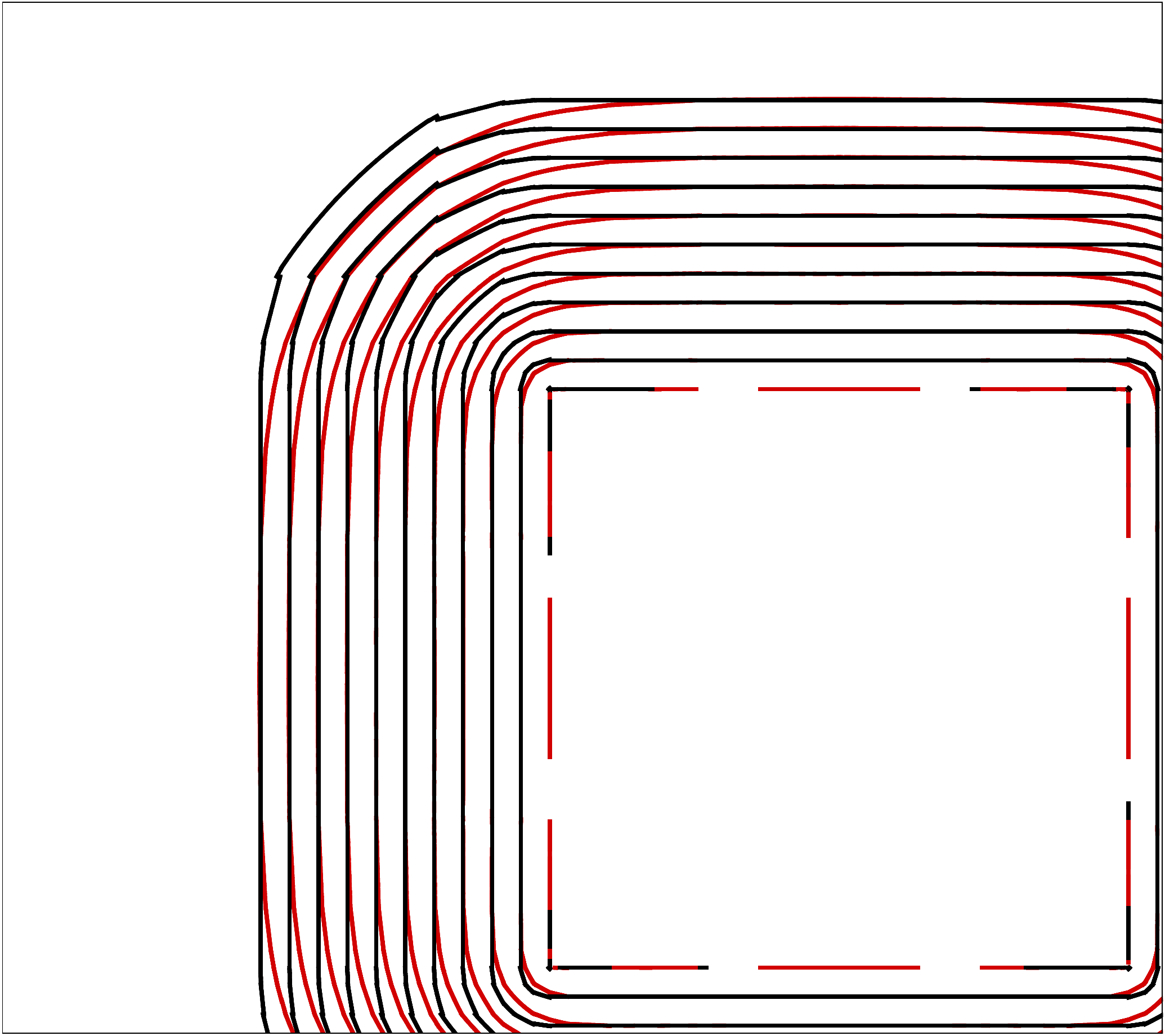}
    \includegraphics[width=0.3\textwidth, trim=80 100 100 10  ,clip]{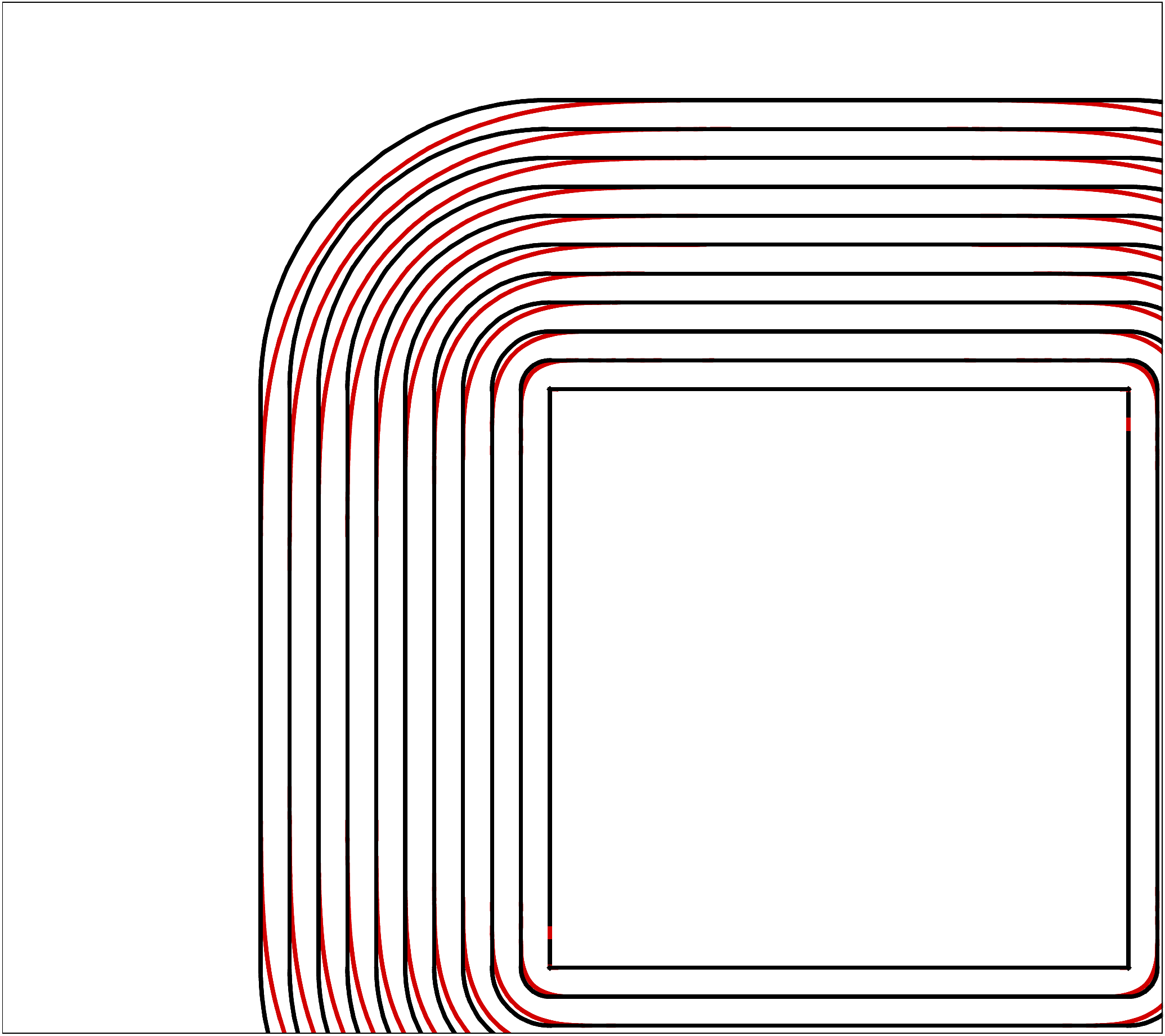}
    \includegraphics[width=0.3\textwidth, trim=80 100 100 10  ,clip]{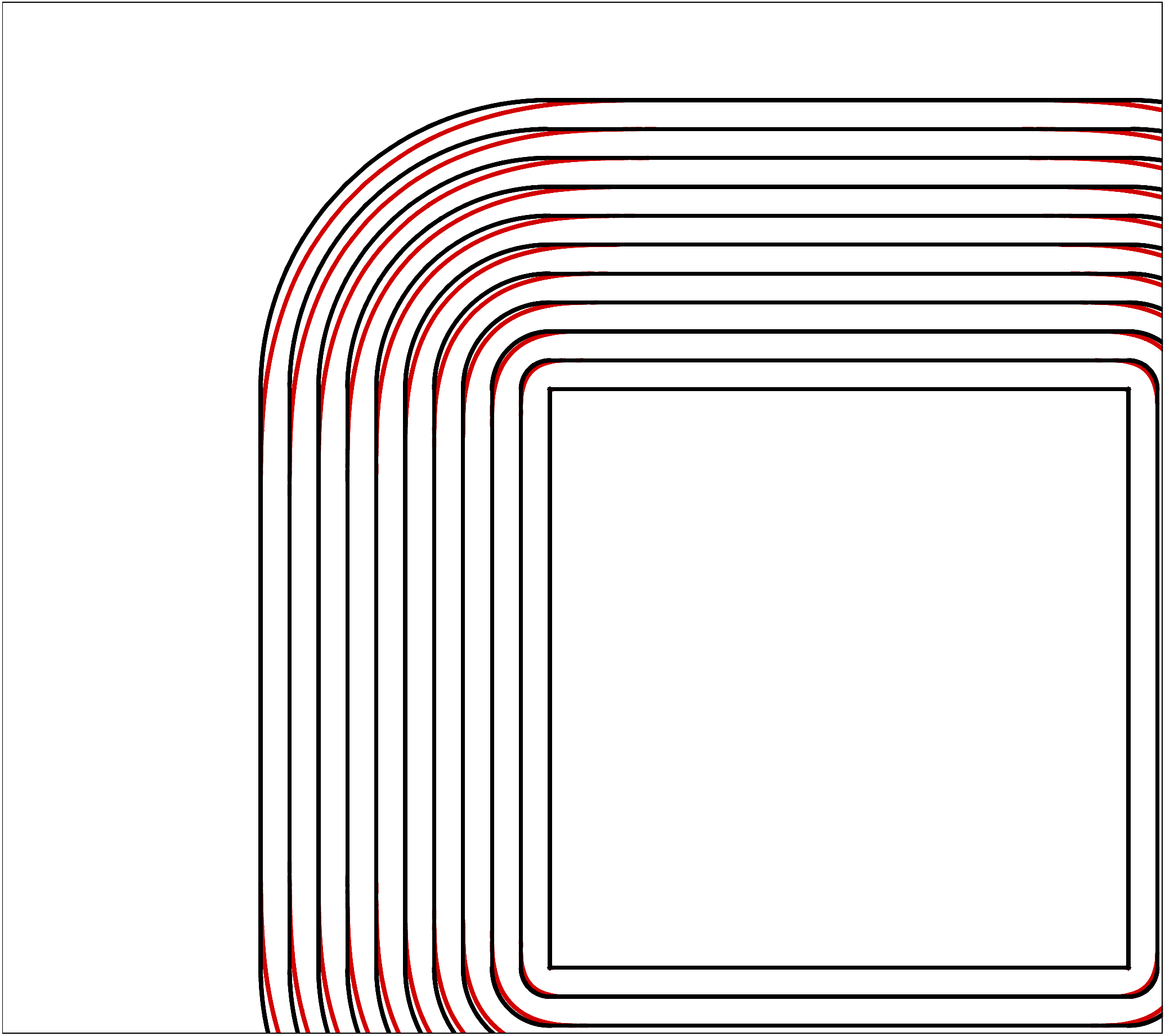}
    \caption{Solution for the square (width $1\times1$) problem for order $N=2,4,8$ with $8$ cells for each side and $8$ cells into the farfield at distance $29$. Red: computed distance function iso-lines [0...1], black: discretization of exact distance, iso-lines [0...1]}
    \label{fig:box}
\end{figure}

In this section the results for the signed distance function around a square are presented. This test illustrates the correct behavior of the proposed method in the case of discontinuous geometries, leading to rarefaction type solutions. We find this to be the most challenging problem. The solution in the rarefaction region is not directly governed by the boundary conditions, as the method does not couple across vertices. In order to have the method converge to the correct answer, the divergence in the source term is scaled by the magnitude of the gradient, $\Vnabla \cdot \vec{q} |\vec{q}|$, wherever the curvature of the gradient field is larger than zero, $\Vnabla \cdot \vec{q} > 0$. Using the curvature is a simple way to distinguish rarefaction and shock regions, by positive or negative curvature respectively, as $s>0$. The exact solution is readily available and we can thus study the convergence of the method in detail. As the geometry has a discontinuity, we cannot expect better than first-order convergence in the asymptotic regime, which indeed is what we find in our numerical experiments. Figure \ref{fig:box} shows the behavior for a coarse H-type mesh around a square, with $8\times8$ elements on each side of the square and into the farfield. The viscosity coefficient is empirically set to $0.3$, and the mesh is linear. The solution clearly opens up in the rarefaction region with increased resolution (order), converging towards the exact solution.

\subsection{Distance from a NACA0012 airfoil}
\begin{figure}
    \centering
        \includegraphics[width=0.48\textwidth, trim=10 10 10 10 ,clip]{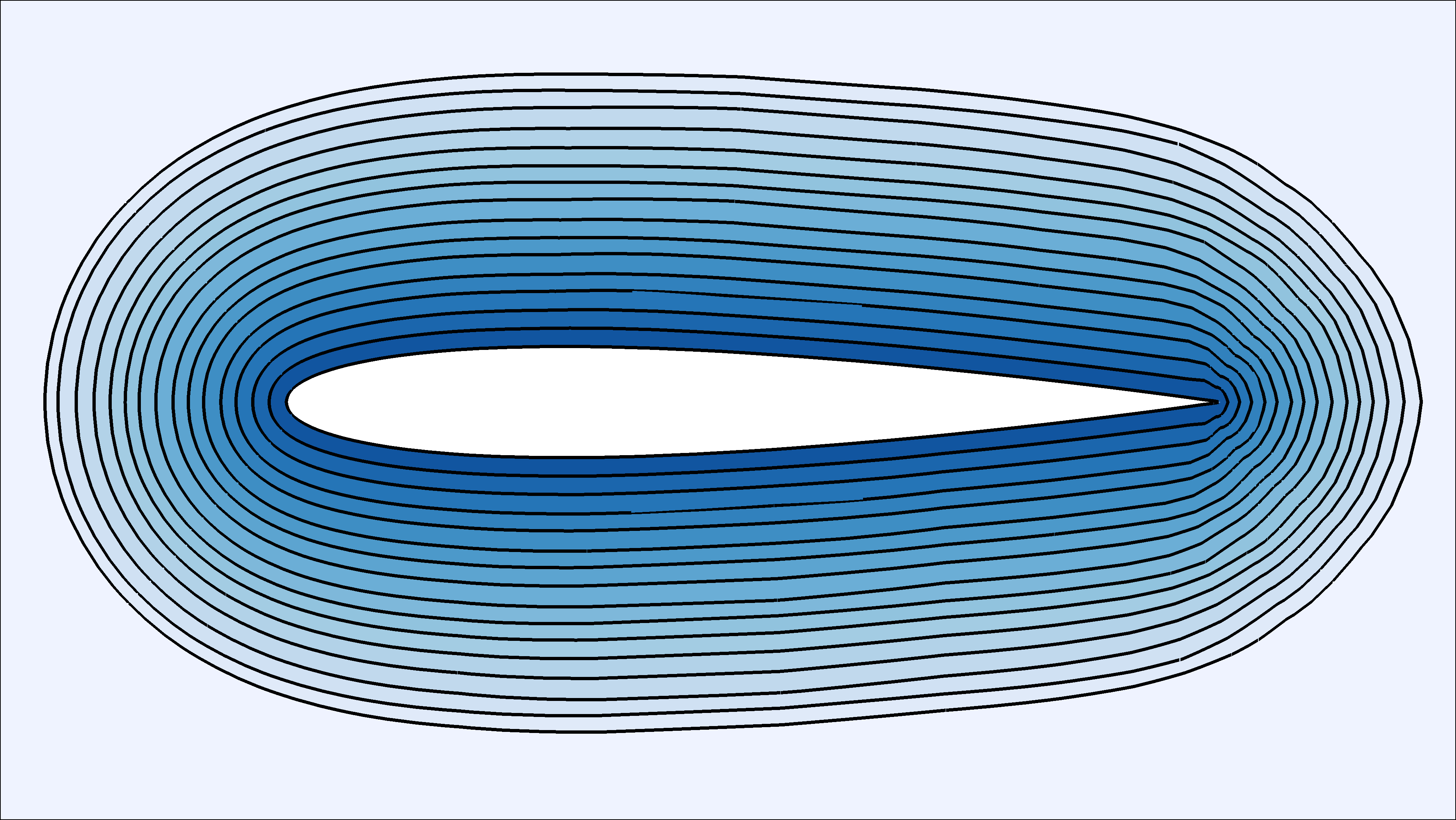}
        \includegraphics[width=0.48\textwidth, trim=10 10 10 10,clip]{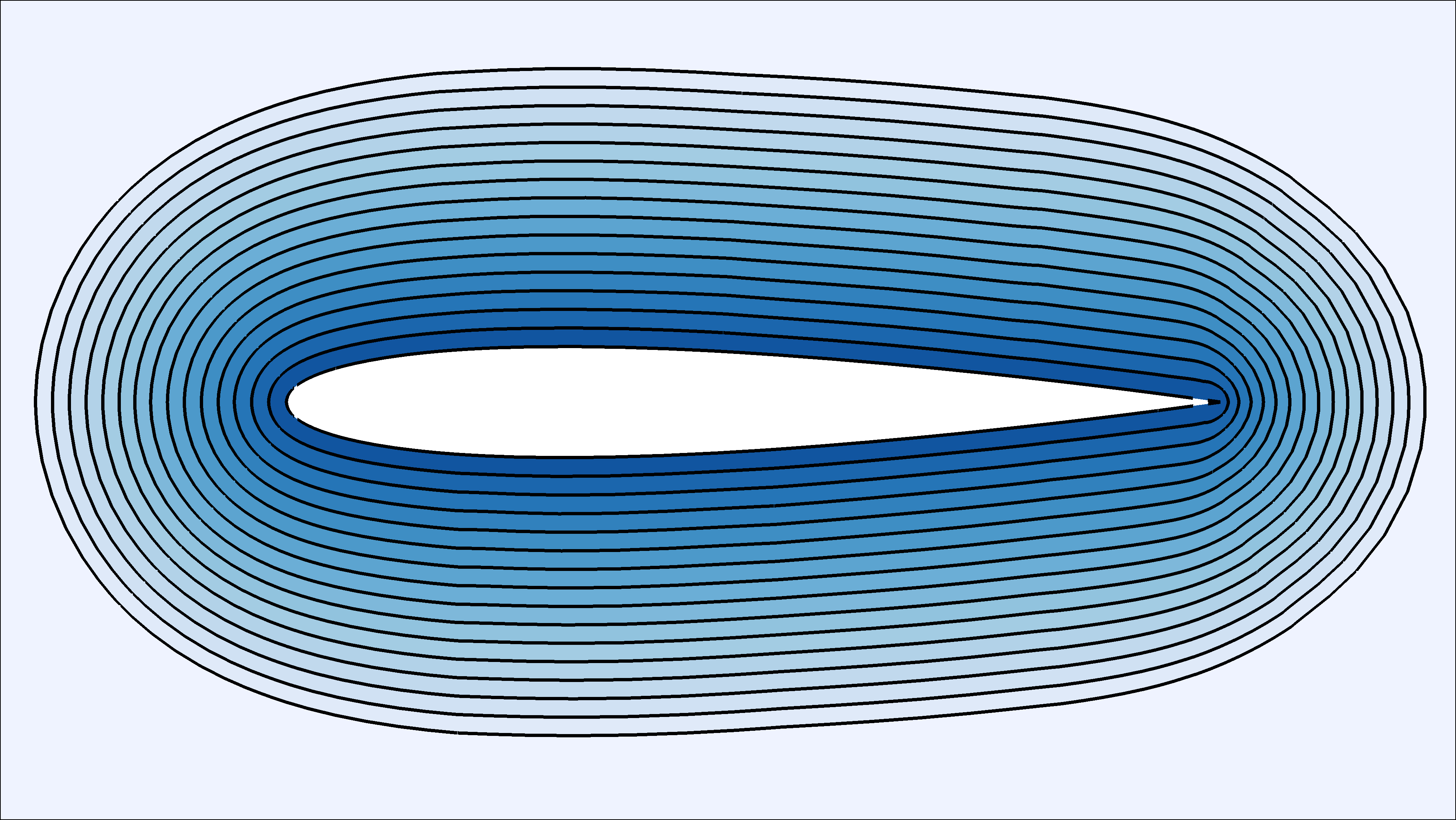}\\
    \includegraphics[width=0.48\textwidth, trim=10 10 10 10 ,clip]{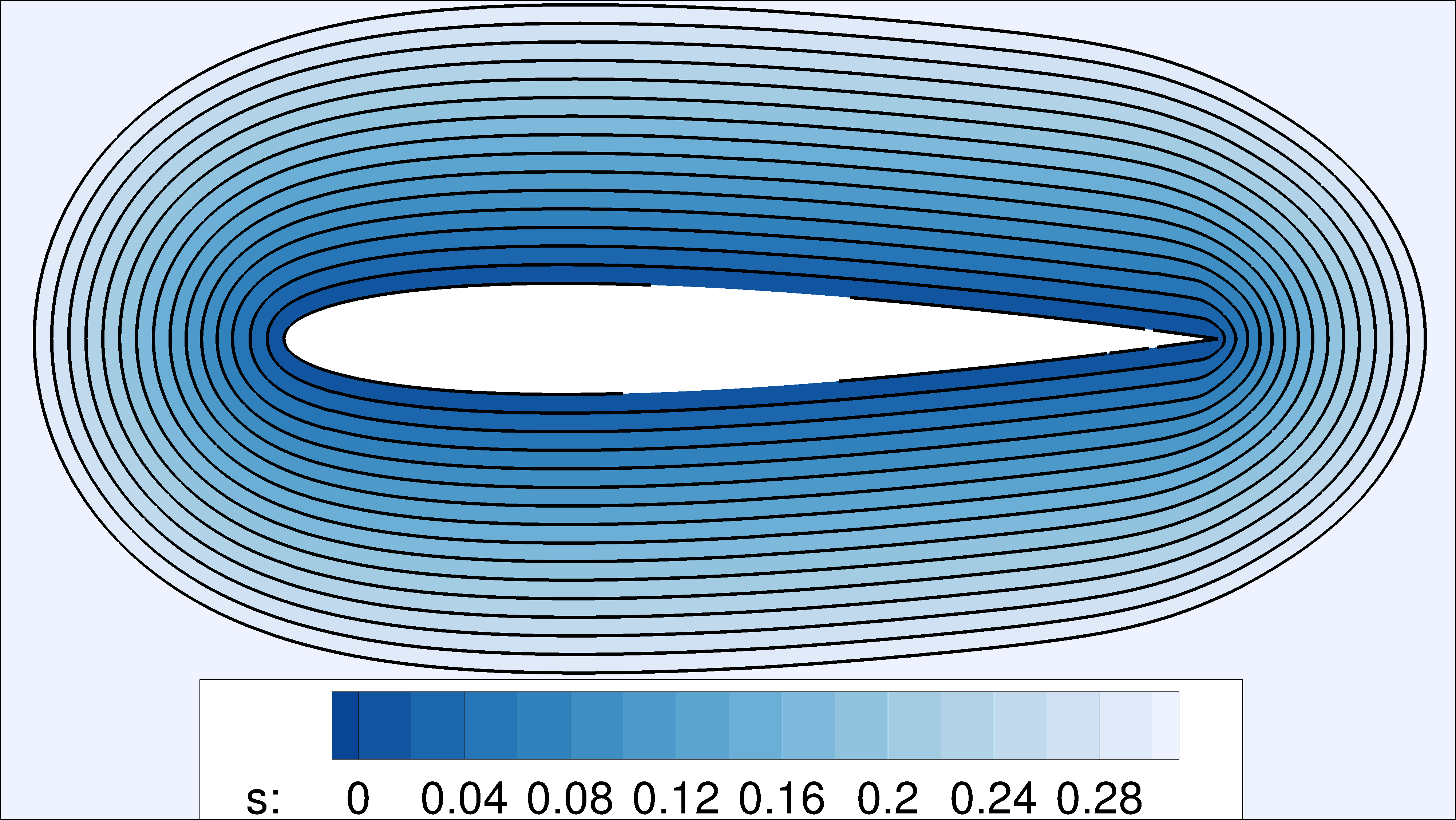}
    \caption{Distance field around a NACA0012 airfoil on $64\times32$ DOF for $N=2,4,8$.}
    \label{fig:naca}
\end{figure}

The distance field around a NACA0012 airfoil is computed using $64\times32$ DOF around the airfoil and in the radial direction. The farfield is circular and located at about $14$ cord length. The computation is done using discretization orders $N=2,4,8$, i.e. the coarsest mesh for $N=8$ is as coarse as $8\times4$ elements. As for the sinusoidal, the artificial viscosity coefficient is set to $c=0.9$. The difficulty of the case is the sharp trailing edge, resulting in a sharp discontinuity of the gradient solution, analogue to a rarefaction. Figure \ref{fig:naca} shows the results of the solution close to the body. For all discretization orders the sharp trailing edge does not lead to an oscillating solution. While all results are qualitatively similar, the high order solutions are considerably smoother.

\subsection{Distance from a 30P30N multi-element airfoil}
\begin{figure}
    \centering
        \includegraphics[width=0.5\textwidth, trim=10 10 10 10 ,clip]{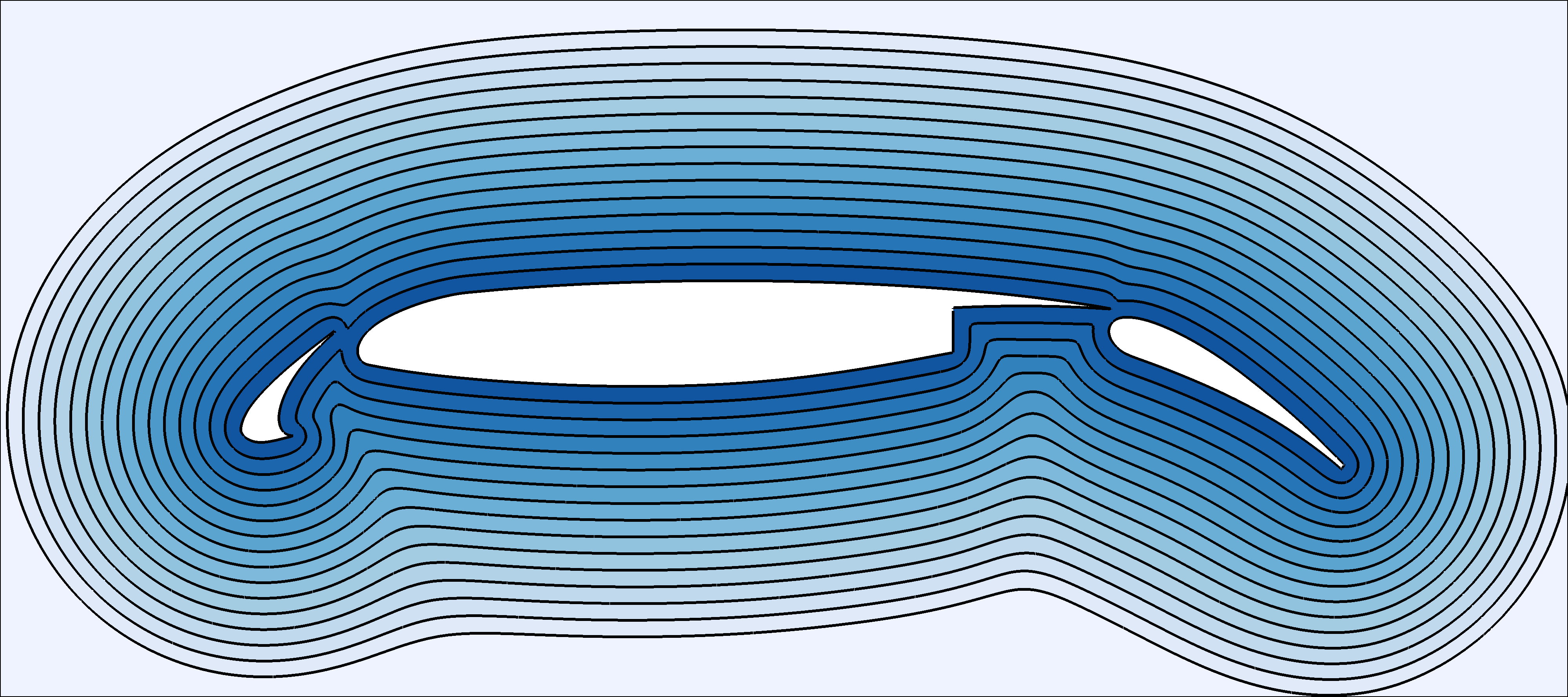}\\
        \includegraphics[width=0.5\textwidth, trim=10 10 10 10,clip]{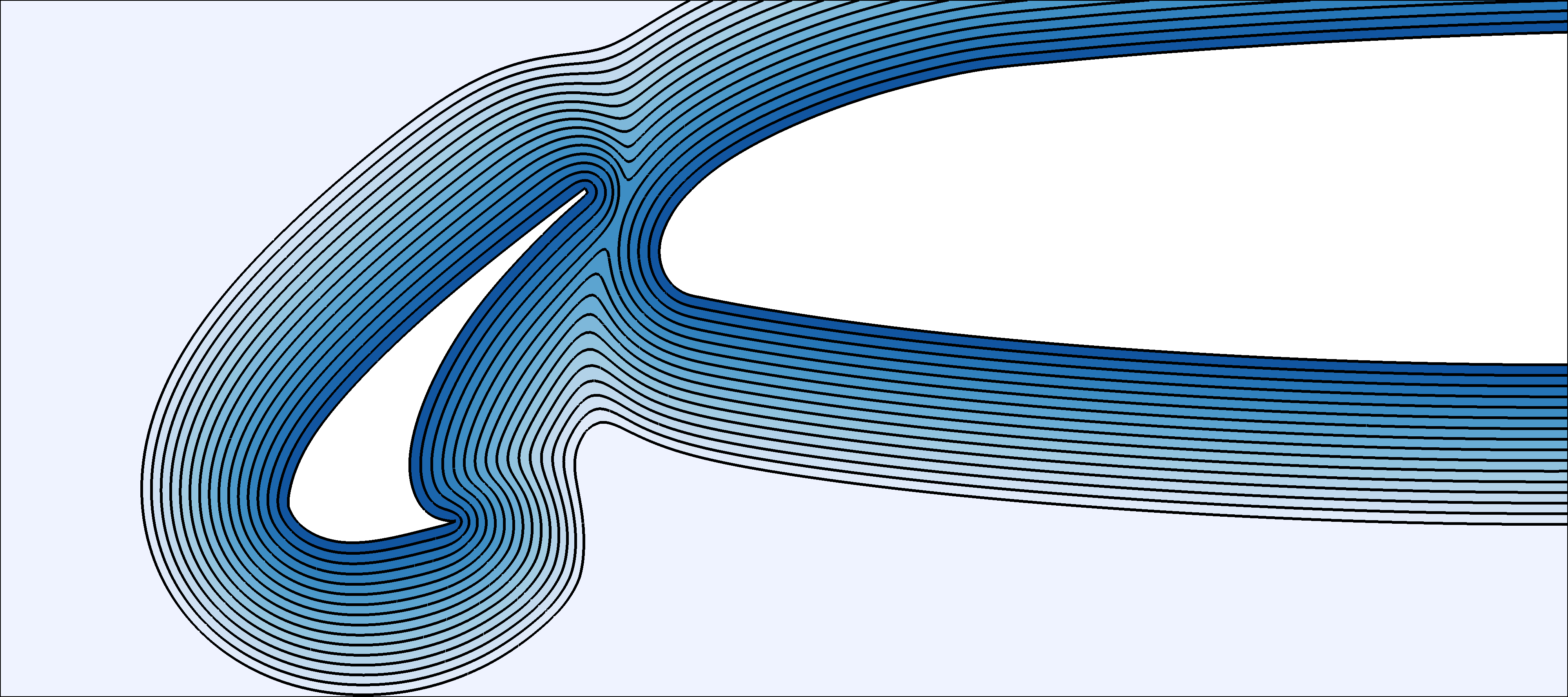}\\
    \includegraphics[width=0.5\textwidth, trim=10 10 10 10 ,clip]{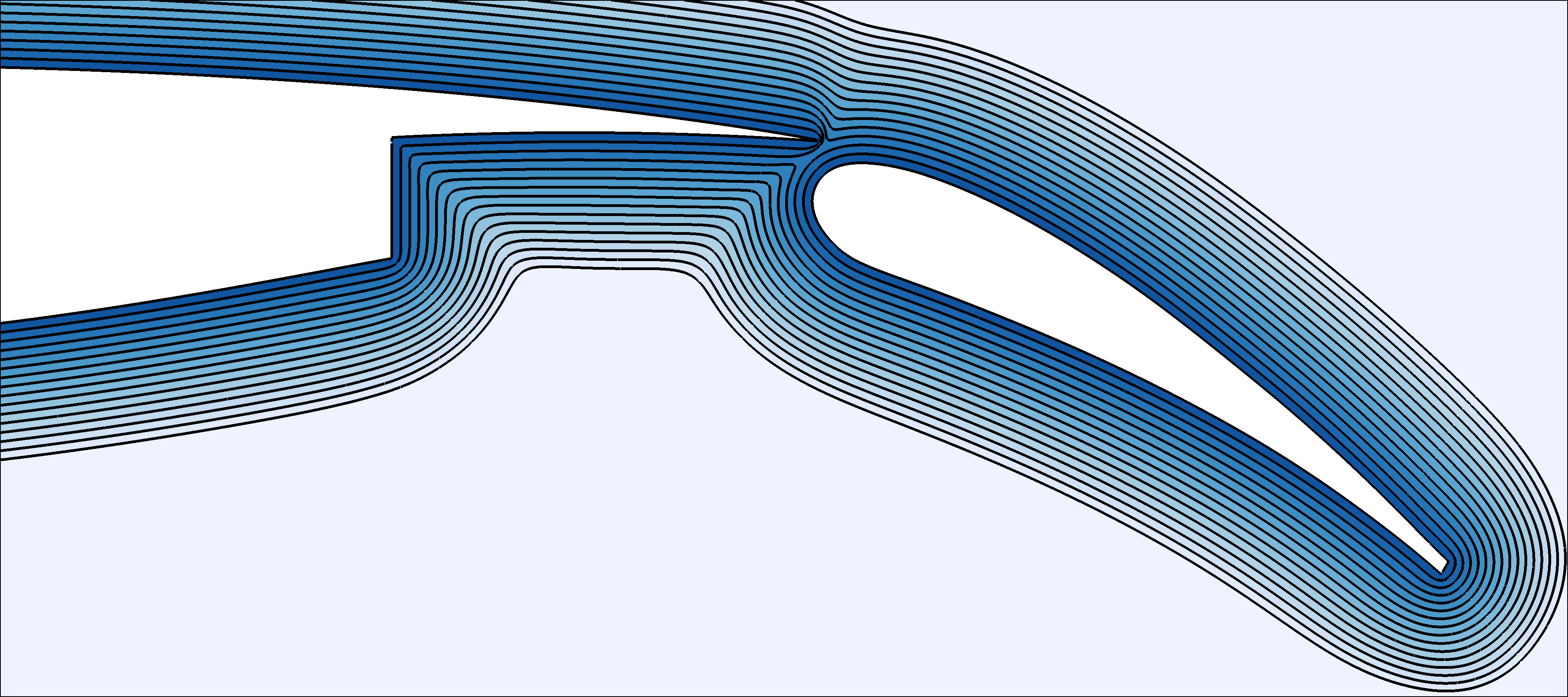}
    \caption{Distance field around a 30P30N, $N=4$, $22,324$ cells. Top: iso-distance lines $[0...5]$, middle and bottom rows: iso-distance lines $[0...1]$.}
    \label{fig:30p30n}
\end{figure}

To demonstrate the capability of the method on a complex geometry, we show results for the NASA 30P30N multi-element airfoil. It is a combination of the above cases, involving concave, convex, smooth and sharp features. A similar airfoil (L1T2) was used by \cite{Hartmann2014} to demonstrate the ability of their continuous-Galerkin eikonal solver to handle complex geometries up to $4^{th}$-order of accuracy. The farfield is circular and located at about $30$ chord length of the main wing. The mesh contains a total of $22,324$ elements and the computation is done for discretization order $N=4$, with a viscosity coefficient $c=0.3$.\footnote{$N=8$ solutions are not presented as the iso-parametric mapping used to generate the mesh did not lead to a quality $9th$-order mesh. Appropriate high-order mesh generation for complex geometry is a topic of ongoing research.}  The resulting distance field is shown in Figure \ref{fig:30p30n}. The distance field is smooth without oscillations, and the high-order solution displays a feature-rich representation of the distance field. Discontinuities at trailing edges are well captured, as are the shocks resulting from the multiple body problem, and the corner inside the flap cove on the main wing. 

\begin{figure}
    \centering
    \includegraphics[width=0.8\textwidth, trim=10 10 10 10 ,clip]{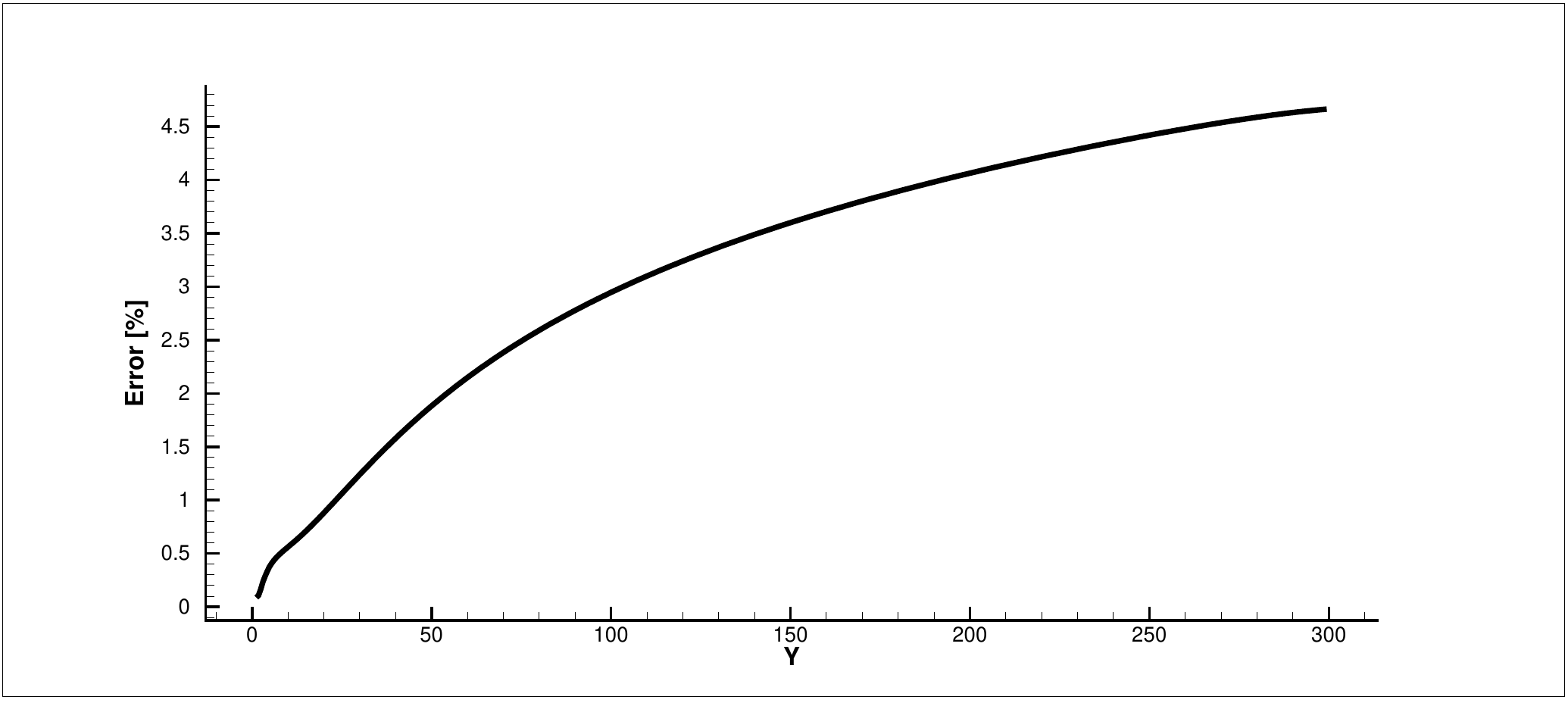}
    \caption{Approximate error of the distance field along a line starting in the center of the main wing.}
    \label{fig:error30P30N}
\end{figure}

To obtain some quantitative measure the distance field is evaluated along a line at the center of the main wing, were the normal vector of the wing aligns with the $Y$ direction of the domain. We compare our solution to the linear $s_{ex} = Y - 1.1145$ in Figure \ref{fig:error30P30N}, plotting the deviation in percent to the linear solution. The approximate error in the near field is well below $1\%$, while it adds up towards $4.7\%$ in the farfield.
%\newpage

\section{Conclusion}
We described a method to compute an approximate distance field around a given obstacle by solving the eikonal equation. Such distance fields are often needed to complement a more relevant physical simulation. Thus, it is desirable to use the same data structure for both, and to have a solution of the same order of accuracy. To obtain a system of equations that is well suited to be solved in our high-order discontinuous-Galerkin framework, we decoupled the gradient solution from the scalar distance field, and only later coupled the two again weakly by adding a source term. For cases where the scalar distance function is discontinuous, leading to shocks, an artificial viscosity term was added. The correct behavior in rarefaction regions is obtained by selectively scaling the driving source term of the system by the magnitude of the gradient.  Formally, this term is equal to one, and hence does not alter the equation.

To demonstrate the capabilities of the proposed method we computed the distance field for the following test cases: a cylinder, showing the method achieves design order of accuracy for smooth solutions; a sinusoidal with a planar wall at the opposite site of the domain, showing the method handles well the appearance of shocks, even when forming complicated structures also in conjunction with convex parts of the solution; a square, showing the method converges towards the entropy solution in the case of a rarefaction; a NACA0012 airfoil with a sharp trailing edge, showing good behavior at very sharp edges even for up to $N=8$ discretizations, and finally, the 30P30N multi-element high-lift configuration with blunt trailing edges, combining the complexity of all the above with a realistically large mesh, showing that the method may be used in actual applications.

\section{Acknowledgement}
David Flad's research was supported by an appointment to the NASA Postdoctoral Program at the NASA Ames Research Center, administered by Universities Space Research Association under contract with NASA.

\newpage
\bibliographystyle{acm}
\bibliography{References}
\end{document}